\documentclass[10pt,a4paper]{amsart}
\usepackage{amssymb}
\usepackage{amsfonts}
\usepackage{amsmath}
\usepackage{epsfig}
\usepackage{indentfirst}
\usepackage[usenames,dvipsnames]{pstricks}
\usepackage{pst-grad}
\usepackage{pst-plot}
\usepackage[margin=1.0in]{geometry}
\usepackage{graphicx}
\usepackage{tikz}
\usepackage{tikz-cd}
\usepackage{float}
\usepackage{tabularx}
\usepackage{mathtools}
\usepackage{bm}
\usepackage{mathrsfs}
\usepackage{lineno}

\newtheorem{theorem}{Theorem}[section]

\newtheorem{corollary}[theorem]{Corollary}

\newtheorem{definition}[theorem]{Definition}

\newtheorem{lemma}[theorem]{Lemma}

\newtheorem{proposition}[theorem]{Proposition}
\newtheorem{remark}[theorem]{Remark}

\newtheorem*{maintheorem}{Theorem A}

\DeclareMathOperator {\Hom}{Hom}

\DeclareMathOperator {\Ind}{Ind}

\newcommand{\Z}{\mathbb{Z}}

\newcommand{\C}{\mathbb{C}}

\newcommand{\B}{\mathcal{B}}

\newcommand{\Sg}{\mathfrak{S}}

\begin{document}
\title[The exact semisimplicity criterion]{The exact semisimplicity criterion for Brauer algebras of type $C$}
\author{Dimpi KM}\address{Indian Institute of Science Education and Research Thiruvananthapuram, Email: dtyagi20@iisertvm.ac.in}

\begin{abstract}
We determine exactly when the Brauer algebra of type $C$, denoted by $\B_r(C,\delta)$, is semisimple over $\C$, where $\delta\neq0$. We first analyze the two admissible cases for homomorphisms from the top layer to the second layer. We construct an auxiliary operator on the space spanned by symmetric Brauer diagrams with a fixed bottom partial diagram and exactly two horizontal edges in each row. Using this operator, we prove that, in each admissible
case, the corresponding quadratic condition is necessary and sufficient for the existence of a nonzero homomorphism. The constant terms of these quadratic expressions depend on sums of contents of horizontal $2$-strips in bipartitions. We then show that a failure of semisimplicity is detected
by cell modules in the second layer. Finally, by relating the content sums for admissible bipartition pairs to the known sums for horizontal $2$-strips in ordinary partitions, we obtain an explicit finite set of exceptional nonzero parameters.
\end{abstract}

\maketitle

\section{Introduction}

The semisimplicity problem for Brauer algebras was studied extensively before the development of cellular algebra theory. In a series of papers, Hanlon and Wales investigated the radicals of Brauer's centralizer algebras and the Wedderburn decompositions
of their semisimple quotients, using the filtration by the number of horizontal edges. They expressed the dimensions of the radicals and of the irreducible representations of the successive quotients in this filtration in terms of the ranks and nullities
of certain matrices constructed from partial matchings. By exploiting the commuting actions of symmetric groups, these large matrices were decomposed into substantially smaller representation-theoretic blocks, whose eigenvalues and determinants control the corresponding radical components. They further developed a tower construction which relates the radicals at consecutive ranks \cite{HW89a,HW89b,HW90,HW94}.

A more uniform framework became available with the introduction of cellular algebras by Graham and Lehrer \cite{GL96}. They showed that the Brauer algebra is cellular and that each cell module carries a canonical associative bilinear form. Consequently, semisimplicity is equivalent to the nondegeneracy of the forms on all cell modules. This formulation organizes the problem through cellular Gram determinants, but it does not by itself determine their exceptional roots.

Doran, Wales, and Hanlon developed a complementary module-theoretic approach using homomorphisms between generic modules, restriction and induction, Frobenius reciprocity, and functors connecting Brauer algebras whose ranks differ by two \cite{DWH99}. Their analysis reduced a possible failure of semisimplicity to content conditions associated with pairs of partitions. Rui later combined these results with cellular theory and proved that, for $\delta\neq0$, the classical Brauer algebra $B_r(\delta)$ over $\C$ is semisimple if and only if all cell modules in the first layer of $B_m(\delta)$ are simple for every $m$ with $2\leq m\leq r$ \cite{Ru05}. Here the
first layer corresponds to diagrams with one horizontal edge in each row. Rui and Si subsequently gave an explicit description of the exceptional set \cite{RS06}, and related methods were later applied to cyclotomic Brauer algebras \cite{RY04, RX07}.

The present work follows this general philosophy for the Brauer algebra of type $C$. The Brauer algebra of type $C$, denoted by $\B_r(C,\delta)$, was introduced by Cohen, Liu and Yu \cite{CLY12}. It can be realized as the subalgebra of the classical Brauer algebra with $2r$ vertices spanned by diagrams which are fixed by reflection about a vertical axis. Bowman proved that this algebra is cellularly stratified \cite{Bo12}. Its cell modules are indexed by a layer and a bipartition, and the top layer is controlled by the hyperoctahedral group $W_r$. These facts make it possible to study semisimplicity through cell modules and homomorphisms between them.

In our earlier work \cite{KCTprevious}, we obtained necessary
conditions for the existence of such homomorphisms. In
particular, if
$$
 \Hom_{\B_r(C,\delta)}
 \bigl(S^{(0,\bm\lambda)},S^{(2,\bm\mu)}\bigr)\neq 0,
$$
then $\delta$ satisfies one of two quadratic
equations whose constant terms are determined by
$\bm\lambda$ and $\bm\mu$. This gave a sufficient condition
for semisimplicity, but not a complete criterion.
Two questions remained. First, are these necessary conditions
also sufficient for the existence of a nonzero homomorphism?
Second, can a deeper layer give an additional obstruction
to semisimplicity?

The present paper resolves both questions. We first study
the second-layer cell modules as $\C W_r$-modules and
describe the two cases in which nonzero homomorphisms of the form
$$
 S^{(0,\bm\lambda)}\longrightarrow S^{(2,\bm\mu)}
$$
can occur. These cases are distinguished by the component
of the bipartition in which a horizontal $2$-strip occurs. We construct an auxiliary operator on the space spanned by symmetric Brauer diagrams with a fixed bottom partial diagram and exactly two horizontal edges in each row. For each bipartition $\bm\mu$ of $r-2$, this operator induces a $\C W_r$-endomorphism $\mathcal G^{\bm\mu}$ on the second-layer
cell module $S^{(2,\bm\mu)}$. Its action on the relevant isotypic components leads to two quadratic expressions:
$$
\delta^2-\delta-2+C_{\bm\lambda/\bm\mu}
\quad\text{and}\quad
\delta^2+\delta-4+C_{\bm\lambda/\bm\mu},
$$
where $C_{\bm\lambda/\bm\mu}$ is twice the sum of the contents of the two added boxes. These expressions correspond, respectively, to a horizontal
$2$-strip in the second and first components of the bipartition.
In each case, the vanishing of the corresponding quadratic
expression is necessary and sufficient for the existence
of a nonzero homomorphism.

We then prove that a failure of semisimplicity is detected
by cell modules in the second layer, see
Theorem~\ref{thm:exact-deciding-layer}.
The isomorphism of homomorphism spaces from our previous
work, recalled in Proposition~\ref{prop:shift-pairing},
allows us to reduce to homomorphisms whose source lies
in the top layer. An argument using induction and
restriction then reduces the obstruction to the second
layer. Together with the preceding characterization of
homomorphisms, this gives a necessary and sufficient
criterion for semisimplicity.

To obtain an explicit criterion on $\delta$, we relate the sums of contents of horizontal $2$-strips in bipartitions to the corresponding known sums for ordinary partitions.
Set
\begin{equation}\label{eq:exact-index-set}
\mathscr S_r= \{s\in\Z:3-r\leq s\leq r-3\}
\cup\{1,3,\ldots,2r-3\},
\end{equation}
where the first set is understood to be empty when $3-r>r-3$. In each of the two cases described above, we show that the possible values of $C_{\bm\lambda/\bm\mu}$ are exactly $2s$ with $s\in\mathscr S_r$. We therefore obtain the following theorem.

\begin{maintheorem}
Let $r\geq 2$ and let $\delta\in\C^{\times}$.  The following conditions are
equivalent.
\begin{enumerate}
\item The algebra $\B_r(C,\delta)$ is semisimple.
\item $ \delta^2-\delta-2+2s\neq 0$ and 
$\delta^2+\delta-4+2s\neq 0$ for every $s\in\mathscr S_r$ 
\end{enumerate}
\end{maintheorem}

Equivalently, the exceptional parameters are exactly the nonzero elements of the finite set
$$
 \left\{\frac{1\pm\sqrt{9-8s}}{2}:s\in\mathscr S_r\right\}
 \cup
 \left\{\frac{-1\pm\sqrt{17-8s}}{2}:s\in\mathscr S_r\right\}.
$$
The formulation in terms of content sums is more convenient for the proofs, while the displayed roots make the exceptional parameters explicit.

The remainder of the paper is organized as follows.
Section~2 recalls the required facts about the hyperoctahedral
group $W_r$, bipartitions, the cell modules of $\B_r(C,\delta)$,
and the results from our previous work. Section~3 determines
the decomposition of the second-layer cell modules as
$\C W_r$-modules and identifies the two admissible cases.
Section~4 describes the two $W_r\times W_{r-2}$-orbits
in the set of symmetric Brauer diagrams with a fixed bottom
partial diagram and exactly two horizontal edges in each row.
It also constructs the auxiliary operator on the span of
these diagrams and the induced operator on each second-layer
cell module $S^{(2,\bm\mu)}$.
Section~5 establishes the precise criterion for the existence
of nonzero homomorphisms from the top layer to the second layer.
Section~6 proves that the second-layer cell modules determine
the semisimplicity of $\B_r(C,\delta)$. Finally, Section~7
determines the possible sums of contents, proves Theorem~A,
and gives the criterion explicitly for $r=2$ and $r=3$.

\section{The Brauer algebra of type $C$ and its cell modules}
\label{sec:background}

Throughout the paper, unless stated otherwise, we work over $\C$. We assume that $\delta\neq0$ and write $\B_r=\B_r(C,\delta)$. All modules are left modules unless stated otherwise. In this section, we recall the notation and results needed for the proof of the main theorem. Further details can be found in \cite{Bo12,CLY12,KCTprevious}.

\subsection{The Brauer algebra of type $C$}

\begin{figure}[H]
\begin{center}
\begin{tikzpicture}[x=0.65cm,y=0.8cm, every node/.style={inner sep=1pt},
    pt/.style={circle,fill=black,inner sep=0pt,minimum size=2pt}]
\node[pt] at (0,0) {};
\node[pt] at (2,0) {};
\node[pt] at (3,0) {};
\node[pt] at (5,0) {};
\node[pt] at (7,0) {};
\node[pt] at (9,0) {};
\node[pt] at (10,0) {};
\node[pt] at (12,0) {};
\node[pt] at (0,-1) {};
\node[pt] at (2,-1) {};
\node[pt] at (3,-1) {};
\node[pt] at (5,-1) {};
\node[pt] at (7,-1) {};
\node[pt] at (9,-1) {};
\node[pt] at (10,-1) {};
\node[pt] at (12,-1) {};
\draw[-] (0,0)-- (0,-1);
\draw[-] (2,0)-- (3,-1);
\draw[-] (3,0)-- (2,-1);
\draw[-] (5,0)-- (5,-1);
\draw[-] (6,0.5)-- (6,-1.5);
\draw[-] (7,0)-- (7,-1);
\draw[-] (9,0)-- (10,-1);
\draw[-] (10,0)-- (9,-1);
\draw[-] (12,0)-- (12,-1);
\node at (-0.75,-0.5) {$s_i=$};
\node at (0,0.5) {\tiny $r$};
\node at (1,-0.5) {\tiny $\cdots$};
\node at (2,0.5) {\tiny $i+1$};
\node at (3,0.5) {\tiny $i$};
\node at (4,-0.5) {\tiny $\cdots$};
\node at (5,0.5) {\tiny $1$};
\node at (7,0.5) {\tiny $-1$};
\node at (8,-0.5) {\tiny $\cdots$};
\node at (8.8,0.5) {\tiny $-i$};
\node at (10.2,0.5) {\tiny $-(i+1)$};
\node at (11,-0.5) {\tiny $\cdots$};
\node at (12,0.5) {\tiny $-r$};
\end{tikzpicture} 
\begin{tikzpicture}[x=0.65cm,y=0.8cm, every node/.style={inner sep=1pt},
    pt/.style={circle,fill=black,inner sep=0pt,minimum size=2pt}]
\node[pt] at (0,0) {};
\node[pt] at (2,0) {};
\node[pt] at (4,0) {};
\node[pt] at (6,0) {};
\node[pt] at (0,-1) {};
\node[pt] at (2,-1) {};
\node[pt] at (4,-1) {};
\node[pt] at (6,-1) {};
\draw[-] (0,0)-- (0,-1);
\draw[-] (2,0)-- (4,-1);
\draw[-] (3,0.5)-- (3,-1.5);
\draw[-] (4,0)-- (2,-1);
\draw[-] (6,0)-- (6,-1);
\node at (-2,-0.5) {$;$};
\node at (-0.75,-0.5) {$t=$};
\node at (0,0.5) {\tiny $r$};
\node at (1,-0.5) {\tiny $\cdots$};
\node at (2,0.5) {\tiny $1$};
\node at (4,0.5) {\tiny $-1$};
\node at (5,-0.5) {\tiny $\cdots$};
\node at (6,0.5) {\tiny $-r$};
\end{tikzpicture}
\begin{tikzpicture}[x=0.65cm,y=0.8cm, every node/.style={inner sep=1pt},
    pt/.style={circle,fill=black,inner sep=0pt,minimum size=2pt}]
\node[pt] at (0,0) {};
\node[pt] at (2,0) {};
\node[pt] at (3,0) {};
\node[pt] at (5,0) {};
\node[pt] at (7,0) {};
\node[pt] at (9,0) {};
\node[pt] at (10,0) {};
\node[pt] at (12,0) {};
\node[pt] at (0,-1) {};
\node[pt] at (2,-1) {};
\node[pt] at (3,-1) {};
\node[pt] at (5,-1) {};
\node[pt] at (7,-1) {};
\node[pt] at (9,-1) {};
\node[pt] at (10,-1) {};
\node[pt] at (12,-1) {};
\draw[-] (0,0)-- (0,-1);
\draw[-] (2,0) to[bend right=40] (3,0);
\draw[-] (2,-1) to[bend left=40] (3,-1);
\draw[-] (5,0)-- (5,-1);
\draw[-] (6,0.5)-- (6,-1.5);
\draw[-] (7,0)-- (7,-1);
\draw[-] (9,0) to[bend right=40] (10,0);
\draw[-] (9,-1) to[bend left=40] (10,-1);
\draw[-] (12,0)-- (12,-1);
\node at (-0.75,-0.5) {$h_i=$};
\node at (0,0.5) {\tiny $r$};
\node at (1,-0.5) {\tiny $\cdots$};
\node at (2,0.5) {\tiny $i+1$};
\node at (3,0.5) {\tiny $i$};
\node at (4,-0.5) {\tiny $\cdots$};
\node at (5,0.5) {\tiny $1$};
\node at (7,0.5) {\tiny $-1$};
\node at (8,-0.5) {\tiny $\cdots$};
\node at (8.8,0.5) {\tiny $-i$};
\node at (10.2,0.5) {\tiny $-(i+1)$};
\node at (11,-0.5) {\tiny $\cdots$};
\node at (12,0.5) {\tiny $-r$};
\end{tikzpicture} 
\begin{tikzpicture}[x=0.65cm,y=0.8cm, every node/.style={inner sep=1pt},
    pt/.style={circle,fill=black,inner sep=0pt,minimum size=2pt}]
\node[pt] at (0,0) {};
\node[pt] at (2,0) {};
\node[pt] at (4,0) {};
\node[pt] at (6,0) {};
\node[pt] at (0,-1) {};
\node[pt] at (2,-1) {};
\node[pt] at (4,-1) {};
\node[pt] at (6,-1) {};
\draw[-] (0,0)-- (0,-1);
\draw[-] (2,0) to[bend right=30] (4,0);
\draw[-] (3,0.5)-- (3,-1.5);
\draw[-] (2,-1) to[bend left=30] (4,-1);
\draw[-] (6,0)-- (6,-1);
\node at (-2,-0.5) {$;$};
\node at (-0.75,-0.5) {$h_0=$};
\node at (0,0.5) {\tiny $r$};
\node at (1,-0.5) {\tiny $\cdots$};
\node at (2,0.5) {\tiny $1$};
\node at (4,0.5) {\tiny $-1$};
\node at (5,-0.5) {\tiny $\cdots$};
\node at (6,0.5) {\tiny $-r$};
\end{tikzpicture}
\caption{Generators of the Brauer algebra of type $C$.}\label{fig:generators-type-C}
\end{center}
\end{figure}
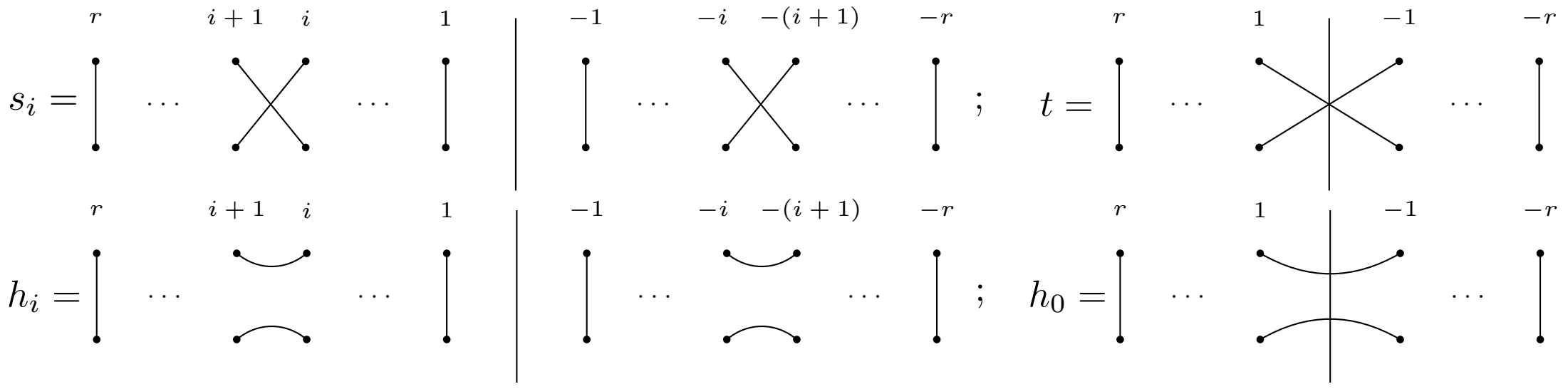

The algebra $\B_r$ is the subalgebra of the classical Brauer algebra $B_{2r}(\delta)$ spanned by the Brauer diagrams fixed under reflection about the vertical axis between the vertices $1$ and $-1$, see \cite{CLY12}. Its generators are represented diagrammatically in Figure~\ref{fig:generators-type-C}.

Let $\mathcal D_r^C$ denote the diagram basis of $\B_r$. The multiplication of basis diagrams is given by concatenation. Each closed loop arising in the middle row is removed and contributes a scalar factor $\delta$. More precisely, if $d_1,d_2\in\mathcal D_r^C$, then $d_1d_2=\delta^{\,j}d_3$, where $j$ is the number of closed loops removed during the concatenation and $d_3\in\mathcal D_r^C$ is the resulting diagram. This multiplication is extended to all of $\B_r$ by $\C$-linearity in each argument.

The diagrams with no horizontal edges form a subalgebra isomorphic to the group algebra of the hyperoctahedral group $W_r$. We therefore identify $\C W_r$ with this subalgebra of $\B_r$. Moreover, the quotient of $\B_r$ by the ideal spanned by the diagrams having at least one horizontal edge is naturally isomorphic to $\C W_r$. 

\subsection{The hyperoctahedral group and bipartitions}

Let $W_r$ be the hyperoctahedral group, realized as the
group of all permutations $w$ of $\{\pm1,\ldots,\pm r\}$
satisfying $w(-i)=-w(i)$ for every $1\leq i\leq r$.
Thus $W_r\cong\Z_2\wr\Sg_r$. For $1\leq i\leq r$,
let $t_i=(i,-i)$, which interchanges $i$ and $-i$
and fixes all other elements. For $1\leq i<j\leq r$,
set
$$
s_{ij}=(i,j)(-i,-j).
$$

The irreducible $\C W_r$-modules are indexed by bipartitions of $r$, see \cite{AK94,Ma99}. A bipartition of $r$ is a pair $\bm\lambda=(\lambda^{(1)},\lambda^{(2)})$ of partitions satisfying $|\lambda^{(1)}|+|\lambda^{(2)}|=r$. We write
$$
 \bm\Lambda_2(r)= \{\bm\lambda=(\lambda^{(1)},\lambda^{(2)}): |\lambda^{(1)}|+|\lambda^{(2)}|=r\}.
$$
The group algebra $\C W_r$ is split semisimple and cellular. For $\bm\lambda\in\bm\Lambda_2(r)$, let $\Delta(\bm\lambda)$ denote the corresponding cell module arising from the cellular construction of Dipper, James and Mathas \cite{DJM98}. Du and Rui proved that the Specht module $S^{\bm\lambda}$ is isomorphic to $\Delta(\bm\lambda)$ \cite{DR01}. Since $\C W_r$ is semisimple, these modules are irreducible. Henceforth, we use the notation $S^{\bm\lambda}$ for the irreducible $\C W_r$-module. 

For a box $p=(i,j)$ in a Young diagram, its content is $ c(p)=j-i$. For $\bm\lambda\in\bm\Lambda_2(r)$, set
\begin{equation}
 C_{\bm\lambda}=2\sum_{p\in[\bm\lambda]}c(p),
 \label{eq:content-scalar}
\end{equation}
where the sum is taken over the boxes in both components of $\bm\lambda$. If $\bm\mu\subseteq\bm\lambda$, set
\begin{equation}
 C_{\bm\lambda/\bm\mu} = C_{\bm\lambda}-C_{\bm\mu} = 2\sum_{p\in[\bm\lambda/\bm\mu]}c(p).
 \label{eq:skew-content}
\end{equation}
Here $\bm\mu\subseteq\bm\lambda$ means
$\mu_i^{(l)}\leq\lambda_i^{(l)}$ for $l=1,2$ and every $i$. The factor $2$ in equation~\eqref{eq:content-scalar} comes from the type $C$ normalization. More precisely, the central element of $\C W_r$ given by
$$
\Omega_r= \sum_{\substack{1\leq i<j\leq r\\ \epsilon=0,1}} t_i^\epsilon t_j^\epsilon s_{ij}
$$
acts on $S^{\bm\lambda}$ by the scalar $C_{\bm\lambda}$; see \cite[Proposition~8.4]{RY04}. We will use this scalar action in the arguments below.

A skew diagram consisting of two boxes is called a \emph{horizontal $2$-strip} if the boxes lie in different columns. The two cases in which such a strip can occur in a bipartition will be identified
in Section~\ref{sec:second-layer-Wr}.

\subsection{Cell modules of the Brauer algebra of type $C$}

For $0\leq k\leq r$, let $P(r,k)$ denote the set of symmetric partial diagrams consisting of a single row of $2r$ vertices with exactly $k$ horizontal edges and $2(r-k)$ free vertices. Each such diagram is invariant under reflection about the vertical axis. Let $V(r,k)$ be the $\C$-vector space spanned by $P(r,k)$. For each $k$, the basis diagrams of $\B_r$ with at least
$k$ horizontal edges in each row span a two-sided ideal $J_k$. Setting $J_{r+1}=0$, we obtain the descending chain of ideals
$$
\B_r=J_0\supseteq J_1\supseteq\cdots\supseteq J_r
\supseteq J_{r+1}=0.
$$
Bowman proved that
\begin{equation}
J_k/J_{k+1}\cong V(r,k)\otimes_{\C}V(r,k)\otimes_{\C}\C W_{r-k}
\label{eq:inflation-layer}
\end{equation}
as an inflation layer \cite{Bo12}. In particular, $\B_r$ is cellularly stratified when $\delta\neq0$.

\begin{figure}[H]
\centering
\begin{tikzpicture}[x=0.6cm,y=0.8cm,
    every node/.style={inner sep=1pt},
    pt/.style={circle,fill=black,inner sep=0pt,minimum size=2pt}
]
\node[pt] at (0,0) {}; 
\node[pt] at (2,0) {};
\node[pt] at (3,0) {}; 
\node[pt] at (4,0) {}; 
\node[pt] at (6,0) {}; 
\node[pt] at (8,0) {}; 
\node[pt] at (10,0) {}; 
\node[pt] at (11,0) {}; 
\node[pt] at (12,0) {};
\node[pt] at (14,0) {}; 
\draw (7,0.5) -- (7,-1.4);
\draw[-] (0,0) to[bend right=14] (14,0);
\draw[-] (2,0) to[bend right=15] (12,0);
\draw[-] (3,0) to[bend right=14] (11,0); 
\node at (-1.5,0) {$D_0=$};
\node at (0,0.5) {\tiny $r$};
\node at (3,0.5) {\tiny $r{-}k{+}1$};
\node at (6,0.5) {\tiny $1$};
\node at (8,0.5) {\tiny $-1$};
\node at (11,0.5) {\tiny $-(r{-}k{+}1)$};
\node at (14,0.5) {\tiny $-r$};
\node at (1,0) {\tiny $\cdots$};
\node at (5,0) {\tiny $\cdots$};
\node at (9,0) {\tiny $\cdots$};
\node at (13,0) {\tiny $\cdots$};

\end{tikzpicture}
\caption{The fixed symmetric partial diagram $D_0$ with $k$ horizontal edges.}
\label{fig:fixed-bottom-partial-diagram}
\end{figure}
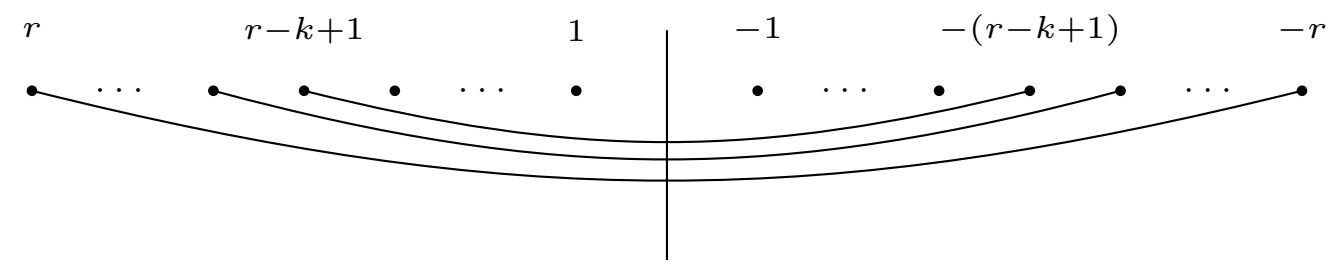

Fix the standard partial diagram $D_0^{(r,k)}\in P(r,k)$ shown in Figure~\ref{fig:fixed-bottom-partial-diagram}. When $r$ and $k$ are clear from the context, we write $D_0$ for $D_0^{(r,k)}$.
As in our previous work \cite{KCTprevious}, we define
$$
 I_r^{k'}(\delta) =V(r,k)\otimes_{\C}\C D_0\otimes_{\C}\C W_{r-k}.
$$
This is a $(\B_r,\C W_{r-k})$-bimodule. For $D\in P(r,k)$, write
$$
 f(D)=D\otimes D_0\otimes1.
$$
Then $\{f(D):D\in P(r,k)\}$ is a basis of $I_r^{k'}(\delta)$ as a
free right $\C W_{r-k}$-module.

The cell modules of $\B_r$ are indexed by
$\Lambda_r=\{(k,\bm\lambda):0\leq k\leq r,\ \bm\lambda\in\bm\Lambda_2(r-k)\}$. We equip $\Lambda_r$ with the partial order used in \cite[Section~2.2]{KCTprevious}: for $(k,\bm\lambda),(l,\bm\mu)\in\Lambda_r$, we define $(k,\bm\lambda)\leq(l,\bm\mu)$ if either $k<l$, or $k=l$ and $\bm\lambda\unlhd\bm\mu$, where $\unlhd$ denotes the dominance order on bipartitions with the convention used in \cite[Section~2.1]{KCTprevious}. For $(k,\bm\lambda)\in\Lambda_r$, the corresponding cell module satisfies the following isomorphism of $\B_r$-modules
\begin{equation}
 S^{(k,\bm\lambda)}\cong I_r^{k'}(\delta)\otimes_{\C W_{r-k}}S^{\bm\lambda}.
 \label{eq:cell-module-realization-new}
\end{equation}
When $k=0$, equation~\eqref{eq:cell-module-realization-new} gives
$S^{(0,\bm\lambda)}\cong S^{\bm\lambda}$. Moreover, $J_1$ annihilates $S^{(0,\bm\lambda)}$, so its action factors through the quotient $\B_r/J_1\cong\C W_r$. Since $S^{\bm\lambda}$ is an irreducible $\C W_r$-module, $S^{(0,\bm\lambda)}$ is a simple $\B_r$-module.

As in \cite[Subsection~2.2]{KCTprevious}, let
$\phi_{(k,\bm\lambda)}$ denote the bilinear form
on $S^{(k,\bm\lambda)}$, and let $G_{(k,\bm\lambda)}$
denote its Gram matrix. By the general theory of cellular algebras
\cite{GL96},
\begin{equation}
\B_r\text{ is semisimple}
\quad\Longleftrightarrow\quad
 \det G_{(k,\bm\lambda)}\neq0
 \quad\text{for every }(k,\bm\lambda)\in\Lambda_r.
 \label{eq:cellular-semisimplicity}
\end{equation}

\subsection{Results from our earlier work \cite{KCTprevious}}

We now recall the results from our previous work \cite{KCTprevious} that will be used in this paper. There, we took a first step towards a semisimplicity criterion by studying homomorphisms between cell modules.

For $\B_r$-modules $M,N$, set
$$
 \langle M,N\rangle_r=\dim_{\C}\Hom_{\B_r}(M,N).
$$

\begin{proposition}[\cite{KCTprevious}, Theorem 4.2]
\label{prop:shift-pairing}
Let $k,k_0\in\mathbb Z_{\geq0}$ with $0\leq k\leq r$, and assume
that $\delta\neq0$. For any bipartitions $\bm{\lambda} \in \bm{\Lambda}_2(r)$, and $\bm{\mu} \in \bm{\Lambda}_2(r-k)$, we have
$$\langle S^{(k_{0}, \bm{\lambda})}, S^{(k+k_{0},\bm{\mu})}\rangle_{r+k_{0}}\neq 0 \quad\Longleftrightarrow\quad \langle S^{(0, \bm{\lambda})}, S^{(k,\bm{\mu})}\rangle_{r}\neq 0. $$
\end{proposition}

If $\bm\alpha$ and $\bm\beta$ are bipartitions, we write $\bm\alpha\to\bm\beta$ if $\bm\beta$ is obtained from $\bm\alpha$ by adding one box to one of its components; equivalently, $\bm\alpha$ is obtained from $\bm\beta$ by removing one box.

\begin{proposition}[\cite{KCTprevious}, Theorem~4.5]
\label{prop:case-split}
Let $1\leq k\leq r$. Suppose that $\bm\lambda\in\bm\Lambda_2(r)$, $\bm\mu\in\bm\Lambda_2(r-k)$, and $\langle S^{(0,\bm\lambda)},S^{(k,\bm\mu)}\rangle_{r}\neq 0$. Then, for every $\bm\lambda^0\to\bm\lambda$, we have
$$
\big\langle S^{(0,\bm\lambda^0)},\,S^{(k,\bm\mu)}\!\downarrow\big\rangle_{r-1}\neq 0,
$$
where $\downarrow$ denotes restriction from $\B_r$ to $\B_{r-1}$. For each such $\bm\lambda^0$, at least one of the following holds; when $k=1$, alternative~(c) is omitted:
\begin{enumerate}
\item $\langle S^{(0,\bm\lambda^0)},\,S^{(k,\bm\mu^1)}\rangle_{r-1}\neq 0$ for some $\bm\mu^1\to\bm\mu$;
\item $\langle S^{(0,\bm\lambda^0)},\,S^{(k-1,\bm\mu)}\rangle_{r-1}\neq 0$;
\item $\langle S^{(0,\bm\lambda^0)},\,S^{(k-2,\bm\mu^0)}\rangle_{r-1}\neq 0$ for some $\bm\mu^0\leftarrow\bm\mu$.
\end{enumerate}
\end{proposition}

Finally, we recall two necessary conditions from our previous
work which motivate the present paper. Before stating them,
we note that
$\langle S^{(0,\bm\lambda)},S^{(2,\bm\mu)}\rangle_r\neq0$
implies $\bm\mu\subseteq\bm\lambda$. Indeed, restricting
a nonzero homomorphism to $\C W_r$ shows that the
irreducible module $S^{\bm\lambda}$ occurs as a constituent
of $S^{(2,\bm\mu)}$. Theorem~3.3 of \cite{KCTprevious},
recalled in
Proposition~\ref{prop:second-layer-Wr-decomposition},
shows that either $\lambda^{(1)}=\mu^{(1)}$ and
$\lambda^{(2)}/\mu^{(2)}$ is a horizontal $2$-strip,
or $\lambda^{(2)}=\mu^{(2)}$ and
$\lambda^{(1)}/\mu^{(1)}$ is a horizontal $2$-strip.
In particular, $\bm\mu\subseteq\bm\lambda$, so the
quantity $C_{\bm\lambda/\bm\mu}$ appearing below is
well defined.

\begin{proposition}\label{prop:known-local-restrictions}
Let $\delta\neq0$.
\begin{enumerate}
\item[(a)] For $\bm\lambda\in\bm\Lambda_2(r)$ and
$\bm\mu\in\bm\Lambda_2(r-1)$, $\langle S^{(0,\bm\lambda)}, S^{(1,\bm\mu)}\rangle_r=0$.
 
\item[(b)] Let $\bm\lambda\in\bm\Lambda_2(r)$ and
$\bm\mu\in\bm\Lambda_2(r-2)$. If $\langle S^{(0,\bm\lambda)}, S^{(2,\bm\mu)}\rangle_r\neq 0$, then either $
 \delta^2-\delta-2+C_{\bm\lambda/\bm\mu}=0$ or $\delta^2+\delta-4+C_{\bm\lambda/\bm\mu}=0$.
\end{enumerate}
\end{proposition}
Part~(a) is \cite[Lemma~5.5]{KCTprevious}, while part~(b) is \cite[Theorem~5.4]{KCTprevious}.

The present paper proves the converse of Proposition~\ref{prop:known-local-restrictions}(b) in the two matched admissible cases. More precisely, we show that, for each of the two horizontal $2$-strip cases, the vanishing of the corresponding quadratic expression gives a nonzero homomorphism from $S^{(0,\bm\lambda)}$ to $S^{(2,\bm\mu)}$. Consequently, the bilinear form $\phi_{(k, \bm\lambda)}$ on $S^{(2,\bm\mu)}$ becomes degenerate. To establish these converses, we first determine the relevant $\C W_r$-module constituents of $S^{(2,\bm\mu)}$. We then use the bimodule $I_r^{2'}(\delta)$ to describe these
constituents diagrammatically and construct the required
auxiliary operator.

\section{The second-layer cell modules as $\C W_r$-modules}
\label{sec:second-layer-Wr}

We determine which irreducible $\C W_r$-modules can occur in a second-layer cell module. This gives a necessary condition for the existence of a nonzero homomorphism from the top layer to the second layer.

Fix $\bm\mu=(\mu^{(1)},\mu^{(2)})\in\bm\Lambda_2(r-2)$. Since
$\C W_r\subseteq\B_r$, we may regard $S^{(2,\bm\mu)}$ as a $\C W_r$-module by restriction. We recall Theorem~3.3 of \cite{KCTprevious}, which describes
the decomposition of $S^{(2,\bm\mu)}$ as a $\C W_r$-module.

\begin{proposition}[\cite{KCTprevious}, Theorem~3.3]
\label{prop:second-layer-Wr-decomposition}
Let $r\geq2$ and $\bm\mu\in\bm\Lambda_2(r-2)$. As a
$\C W_r$-module,
$$
S^{(2,\bm\mu)}
\cong\bigoplus_{\bm\lambda\in\bm\Lambda_2(r)}a_{\bm\mu}^{\bm\lambda}S^{\bm\lambda},
$$
where
$a_{\bm\mu}^{\bm\lambda}=2c_{\bm\mu}^{\bm\lambda}+d_{\bm\mu}^{\bm\lambda}$. Here
$$
c_{\bm\mu}^{\bm\lambda}=\begin{cases}
1, & \lambda^{(2)}=\mu^{(2)}
\text{ and } \lambda^{(1)}/\mu^{(1)}
\text{ is a horizontal $2$-strip},\\
0, & \text{otherwise},
\end{cases}
$$
and
$$
d_{\bm\mu}^{\bm\lambda} = \begin{cases}
1, & \lambda^{(1)}=\mu^{(1)}
\text{ and } \lambda^{(2)}/\mu^{(2)}
\text{ is a horizontal $2$-strip},\\
0, & \text{otherwise}.
\end{cases}
$$
\end{proposition}

The conditions $c_{\bm\mu}^{\bm\lambda}=1$ and
$d_{\bm\mu}^{\bm\lambda}=1$ cannot hold simultaneously.
Indeed, the first implies
$|\lambda^{(1)}|-|\mu^{(1)}|=2$, whereas the second
requires $\lambda^{(1)}=\mu^{(1)}$. Consequently,
$a_{\bm\mu}^{\bm\lambda}
=2c_{\bm\mu}^{\bm\lambda}+d_{\bm\mu}^{\bm\lambda}
\in\{0,1,2\}$. For a $\C W_r$-module $M$, let $[M:S^{\bm\lambda}]_{\C W_r}$ denote the multiplicity of $S^{\bm\lambda}$ in $M$. Proposition~\ref{prop:second-layer-Wr-decomposition}
gives
$$
\bigl[S^{(2,\bm\mu)}:S^{\bm\lambda}\bigr]_{\C W_r}=
\begin{cases}
2, & c_{\bm\mu}^{\bm\lambda}=1,\\
1, & d_{\bm\mu}^{\bm\lambda}=1,\\
0, & c_{\bm\mu}^{\bm\lambda}=d_{\bm\mu}^{\bm\lambda}=0.
\end{cases}
$$

For $\bm\lambda\in\bm\Lambda_2(r)$, the restriction of $S^{(0,\bm\lambda)}$ to $\C W_r$ is the irreducible module
$S^{\bm\lambda}$. Moreover, every $\B_r$-homomorphism remains a $\C W_r$-homomorphism after restriction. Therefore, restriction gives an
injective linear map
$$
\Hom_{\B_r}\bigl(S^{(0,\bm\lambda)},S^{(2,\bm\mu)}\bigr)\longrightarrow
\Hom_{\C W_r}\bigl(S^{\bm\lambda},S^{(2,\bm\mu)}\bigr).
$$
Since $\C W_r$ is semisimple, the dimension of the space on the right is the multiplicity of $S^{\bm\lambda}$ in $S^{(2,\bm\mu)}$. Hence
$$
\langle S^{(0,\bm\lambda)},S^{(2,\bm\mu)}
\rangle_r \leq a_{\bm\mu}^{\bm\lambda}.
$$
This necessary condition leads to the following definition.

\begin{definition}
\label{def:mu-admissible}
Let $\bm\mu\in\bm\Lambda_2(r-2)$. A bipartition $\bm\lambda\in\bm\Lambda_2(r)$ is called
\emph{$\bm\mu$-admissible} if $a_{\bm\mu}^{\bm\lambda}>0$. We denote the set of all $\bm\mu$-admissible bipartitions by
$$
\mathcal A(\bm\mu)=\left\{
\bm\lambda\in\bm\Lambda_2(r)\mathrel{\Big|}a_{\bm\mu}^{\bm\lambda}>0\right\}.
$$
We further define
$$
\mathcal A^-(\bm\mu)=\left\{
\bm\lambda\in\bm\Lambda_2(r)\mathrel{\Big|}\lambda^{(1)}=\mu^{(1)}\text{ and }
\lambda^{(2)}/\mu^{(2)}\text{ is a horizontal $2$-strip}\right\}
$$
and
$$
\mathcal A^+(\bm\mu)=\left\{
\bm\lambda\in\bm\Lambda_2(r)\mathrel{\Big|}\lambda^{(2)}=\mu^{(2)}\text{ and }
\lambda^{(1)}/\mu^{(1)}\text{ is a horizontal $2$-strip}\right\}.
$$
\end{definition}

The two sets are disjoint, and $\mathcal A(\bm\mu)=\mathcal A^-(\bm\mu)\mathbin{\dot\cup}\mathcal A^+(\bm\mu)
$. Moreover, $\bm\lambda\in\mathcal A^-(\bm\mu)$ if and only if $d_{\bm\mu}^{\bm\lambda}=1$, whereas $\bm\lambda\in\mathcal A^+(\bm\mu)$ if and only if $c_{\bm\mu}^{\bm\lambda}=1$. It follows from the preceding inequality that $\langle S^{(0,\bm\lambda)},S^{(2,\bm\mu)}\rangle_r\neq0$ only if
$\bm\lambda\in\mathcal A(\bm\mu)$. Admissibility alone does not guarantee the existence of a nonzero $\B_r$-homomorphism. It remains to determine whether the $S^{\bm\lambda}$-isotypic component of $S^{(2,\bm\mu)}$ contains a $\B_r$-submodule isomorphic
to $S^{(0,\bm\lambda)}$. To study this question, we use the bimodule $I_r^{2'}(\delta)$ and introduce an an auxiliary operator, which will be used to determine when such $\B_r$-submodules exist.

\section{The fixed-bottom space and the auxiliary operator}
\label{sec:fixed-bottom-operator}
The decomposition of $S^{(2,\bm\mu)}$ as a $\C W_r$-module identifies its irreducible constituents, but it does not locate the corresponding submodules inside the cell module. To study these submodules, we use the bimodule $I_r^{2'}(\delta)$ introduced in Section~\ref{sec:background}, which provides a diagrammatic description of the cell module.

Fix $\bm\mu\in\bm\Lambda_2(r-2)$, and let
$D_0=D_0^{(r,2)}$ be the fixed partial diagram shown in
Figure~\ref{fig:fixed-bottom-partial-diagram} with $k=2$. By the construction in Section~\ref{sec:background},
$$
I_r^{2'}=\bigoplus_{D\in P(r,2)}f(D)R
$$
is a $(\B_r,\C W_{r-2})$-bimodule.
We refer to $I_r^{2'}(\delta)$ as the \emph{fixed-bottom space},
since every diagram in its natural basis has bottom partial
diagram $D_0$. By the isomorphism
\eqref{eq:cell-module-realization-new}, we have
$S^{(2,\bm\mu)}
\cong I_r^{2'}(\delta)\otimes_{\C W_{r-2}}S^{\bm\mu}$.
If $\{v_1,\ldots,v_d\}$ is a basis of $S^{\bm\mu}$, then
$$
\left\{
f(D)\otimes_{\C W_{r-2}}v_j
\mid D\in P(r,2),\ 1\leq j\leq d
\right\}
$$
forms a $\C$-basis of $S^{(2,\bm\mu)}$.

For $x\in W_r$ and $D\in P(r,2)$, there exist unique elements
$xD\in P(r,2)$ and $\sigma(x,D)\in W_{r-2}$ such that
$$
xf(D)=f(xD)\sigma(x,D).
$$
Here, $xD$ is obtained by applying the signed permutation $x$
to the vertices of $D$, while $\sigma(x,D)$ records the
induced permutation and sign changes on the remaining
vertical edges. Consequently, the action of $x$ on the
cell module is given by
$$
x\bigl(f(D)\otimes_{\C W_{r-2}}v\bigr)
=f(xD)\otimes_{\C W_{r-2}}\sigma(x,D)v.
$$

\subsection{The two diagram orbits}

The group $W_r\times W_{r-2}$ acts on the diagram basis of
$I_r^{2'}(\delta)$ by $(x,w)\cdot y=xyw^{-1}$. As in \cite[Section~3]{KCTprevious}, this basis decomposes
into two orbits, denoted by $D_{(2,0)}$ and $D_{(0,1)}$.
For $(a,b)\in\{(2,0),(0,1)\}$, define
$$
P_{(a,b)}
=\left\{
D\in P(r,2)\mid f(D)\in D_{(a,b)}
\right\}.
$$
These are the two $W_r$-orbits in $P(r,2)$, and
$P(r,2)=P_{(2,0)}\mathbin{\dot\cup}P_{(0,1)}$. In $P_{(2,0)}$, each horizontal edge is fixed by reflection,
so the two edges have the form $\{i,-i\}$ and $\{j,-j\}$.
In $P_{(0,1)}$, reflection interchanges the two horizontal
edges. They have either the form $\{i,j\}$ and $\{-i,-j\}$,
or the form $\{i,-j\}$ and $\{-i,j\}$, where
$1\leq i<j\leq r$.

Define
$$
I_{20}:=\C D_{(2,0)}
=\bigoplus_{D\in P_{(2,0)}}f(D)\C W_{r-2},
\qquad
I_{01}:=\C D_{(0,1)}
=\bigoplus_{D\in P_{(0,1)}}f(D)\C W_{r-2},
$$
where $\C D_{(2,0)}$ and $\C D_{(0,1)}$ denote the complex
linear spans of the respective orbits. Both spaces are stable
under the left $\C W_r$-action and the right
$\C W_{r-2}$-action. Hence
$$
I_r^{2'}(\delta)=I_{20}\oplus I_{01}
$$
as $(\C W_r,\C W_{r-2})$-bimodules. Tensoring with
$S^{\bm\mu}$ gives
$$
S^{(2,\bm\mu)}
\cong
\bigl(I_{20}\otimes_{\C W_{r-2}}S^{\bm\mu}\bigr)
\oplus
\bigl(I_{01}\otimes_{\C W_{r-2}}S^{\bm\mu}\bigr)
$$
as $\C W_r$-modules. These two summands need not be
$\B_r$-submodules.

We now recall how the orbit spaces separate the irreducible constituents. From \cite[Section~3]{KCTprevious}, we have the following
decompositions into irreducible $\C W_r$-modules
$$
I_{20}\otimes_R S^{\bm\mu}\cong
\bigoplus_{\bm\lambda\in\bm\Lambda_2(r)}
c_{\bm\mu}^{\bm\lambda}S^{\bm\lambda}
$$
and
$$
I_{01}\otimes_R S^{\bm\mu}\cong
\bigoplus_{\bm\lambda\in\bm\Lambda_2(r)}
\bigl(c_{\bm\mu}^{\bm\lambda}
+d_{\bm\mu}^{\bm\lambda}\bigr)S^{\bm\lambda}.
$$
These decompositions locate the admissible constituents. Suppose if $\bm\lambda\in\mathcal A^-(\bm\mu)$,
the unique copy of $S^{\bm\lambda}$ occurs in the summand $I_{01}\otimes_{\C W_{r-2}}S^{\bm\mu}$.
Denote this submodule by
$$
M_{01}^{\bm\lambda,\bm\mu}
\subseteq I_{01}\otimes_{\C W_{r-2}} S^{\bm\mu}.
$$
If $\bm\lambda\in\mathcal A^+(\bm\mu)$, then $S^{\bm\lambda}$ occurs once in each of the two orbit summands. Denote these submodules by
$$
M_{20}^{\bm\lambda,\bm\mu}
\subseteq I_{20}\otimes_{\C W_{r-2}} S^{\bm\mu},
\qquad
M_{01}^{\bm\lambda,\bm\mu}
\subseteq I_{01}\otimes_{\C W_{r-2}} S^{\bm\mu}.
$$
In this case, their direct sum is the full $\bm\lambda$-isotypic component of $S^{(2,\bm\mu)}$.

For a basis diagram $y\in D_{(0,1)}$, define $\tau(y)$
by interchanging the two possible configurations of its
top horizontal edges:
$$
\{i,j\},\{-i,-j\}
\quad\longleftrightarrow\quad
\{i,-j\},\{-i,j\},
\qquad 1\leq i<j\leq r,
$$
and leaving all other edges unchanged. Extend this assignment
$\C$-linearly to a map
$$
\tau:I_{01}\longrightarrow I_{01}.
$$

\begin{lemma}\label{lem:exact-tau-action}
The map $\tau$ is an involution of the
$(\C W_r,\C W_{r-2})$-bimodule $I_{01}$. It induces a
$\C W_r$-linear involution on
$I_{01}\otimes_{\C W_{r-2}}S^{\bm\mu}$, which we also
denote by $\tau$. Moreover,
$$
\tau\big|_{M_{01}^{\bm\lambda,\bm\mu}}
=
\begin{cases}
-\operatorname{Id},
& \bm\lambda\in\mathcal A^-(\bm\mu),\\
\operatorname{Id},
& \bm\lambda\in\mathcal A^+(\bm\mu).
\end{cases}
$$
\end{lemma}

\begin{proof}
By definition, $\tau^2=\operatorname{Id}$. A signed
permutation maps the two configurations with endpoints
$\{\pm i,\pm j\}$ to the two configurations on the image
of this set. It also induces the same permutation and sign
changes on the remaining vertical edges in both configurations.
Thus
$$
\tau(xy)=x\tau(y)
\qquad(x\in W_r,\ y\in D_{(0,1)}).
$$
The right action of $W_{r-2}$ changes only the lower endpoints
of the vertical edges and therefore commutes with $\tau$.
Hence $\tau$ is a
$(\C W_r,\C W_{r-2})$-bimodule involution.

In the orbit decomposition used in
\cite[Section~3]{KCTprevious}, the two local configurations
give the permutation representation
$$
\Ind_K^{W_2}1\cong1\oplus\chi,
$$
where $K=\langle s_{12},t_1t_2\rangle$ has index two in
$W_2$, and $\chi$ is the nontrivial one-dimensional
representation of $W_2/K$. The operation defining $\tau$
interchanges the two cosets of $K$. It therefore acts as
$+1$ on the trivial summand and as $-1$ on $\chi$.
The trivial summand accounts for the constituents whose
multiplicities are given by $c_{\bm\mu}^{\bm\lambda}$,
whereas $\chi$ accounts for those whose multiplicities
are given by $d_{\bm\mu}^{\bm\lambda}$.
The stated action follows.
\end{proof}

\subsection{The auxiliary operator}

Let $*$ denote the diagram anti-involution obtained by
reflection across a horizontal axis. On the standard basis arising from the inflation construction, it is given by
$$
\bigl(D_1\otimes D_2\otimes w\bigr)^*
=D_2\otimes D_1\otimes w^{-1}.
$$
In particular,
$$
f(D)^*=D_0\otimes D\otimes 1.
$$
For $D,E\in P(r,2)$, define $\phi(D,E)\in\C W_{r-2}$ by
\begin{equation}\label{eq:exact-phi-definition}
f(D)^*f(E)=D_0\otimes D_0\otimes\phi(D,E)
\end{equation}
in the quotient $J_2/J_3$. If the concatenation produces more than two horizontal edges, its image in this quotient is zero, and we set
$\phi(D,E)=0$. Otherwise,
$$
\phi(D,E)=\delta^{\,t(D,E)}w(D,E)
$$
where $t(D,E)\in\{0,1,2\}$ is the number of closed loops removed during concatenation and $w(D,E)\in W_{r-2}$ is the signed permutation determined by the remaining vertical edges. We also use $*$ to denote the $\C$-linear anti-involution
of $\C W_{r-2}$ defined by $w^*=w^{-1}$ for
$w\in W_{r-2}$. Applying $*$ to
\eqref{eq:exact-phi-definition} gives
$$
\phi(D,E)^*=\phi(E,D).
$$
Moreover, $\phi(D,D)=\delta^2 1_{\C W_{r-2}}$, and if $D\neq E$, then $\phi(D,E)\in\{0,w,\delta w\}$ for some $w\in W_{r-2}$.

Let $P(r,2)=\{D_1,\ldots,D_N\}$. The coefficients $\phi(D_i,D_j)$ form the $\C W_{r-2}$-valued matrix
$$
\mathcal G_r(\delta) := \bigl(\phi(D_i,D_j)\bigr)_{1\leq i,j\leq N}.
$$
This matrix defines a right $\C W_{r-2}$-linear operator
$$
\mathcal G_r(\delta):
I_r^{2'}(\delta)\longrightarrow I_r^{2'}(\delta)
$$
by
$$
\mathcal G_r(\delta)\bigl(f(E)\bigr) = \sum_{D\in P(r,2)}f(D)\phi(D,E)
$$
and, more generally,
$$
\mathcal G_r(\delta)\bigl(f(E)a\bigr)
=
\sum_{D\in P(r,2)}f(D)\phi(D,E)a,
\qquad a\in \C W_{r-2}.
$$
This operator is well defined because every element of $I_r^{2'}(\delta)$ has a unique expansion with respect to the free right $\C W_{r-2}$-basis $\{f(D):D\in P(r,2)\}$.

\begin{proposition}
\label{prop:auxiliary-operator-bimodule-map}
The operator $\mathcal G_r(\delta)$ is a $(\C W_r,\C W_{r-2})$-bimodule endomorphism of $I_r^{2'}(\delta)$.
\end{proposition}

\begin{proof}
By definition, $\mathcal G_r(\delta)$ is right $\C W_{r-2}$-linear. It remains to show that it is left $\C W_r$-linear. Let $x\in W_r$ and
$D,E\in P(r,2)$. Recall that
$$
xf(D)=f(xD)\sigma(x,D), \qquad xf(E)=f(xE)\sigma(x,E).
$$
Since $x^*=x^{-1}$, we have
$$
\bigl(xf(D)\bigr)^*\bigl(xf(E)\bigr)=f(D)^*f(E).
$$
Using the definition of $\phi$, we obtain
$$
\sigma(x,D)^{-1}\phi(xD,xE)\sigma(x,E)=\phi(D,E),
$$
and hence
\begin{equation}\label{eq:exact-phi-equivariance}
\phi(xD,xE)\sigma(x,E) = \sigma(x,D)\phi(D,E).
\end{equation}
We now verify that $\mathcal G_r(\delta)$ commutes with the left action of $W_r$. By definition,
$$
\begin{aligned}
\mathcal G_r(\delta)\bigl(xf(E)\bigr)
&=\sum_{D'\in P(r,2)}
f(D')\phi(D',xE)\sigma(x,E)\\
&=\sum_{D\in P(r,2)}
f(xD)\phi(xD,xE)\sigma(x,E)\\
&=\sum_{D\in P(r,2)}
f(xD)\sigma(x,D)\phi(D,E)\\
&=x\mathcal G_r(\delta)\bigl(f(E)\bigr),
\end{aligned}
$$
Where the second equality follows by reindexing the sum using $D'=xD$, and the third follows from
\eqref{eq:exact-phi-equivariance}. Thus
$\mathcal G_r(\delta)$ is also left $\C W_r$-linear and therefore is a $(\C W_r,\C W_{r-2})$-bimodule endomorphism.
\end{proof}

Since $\mathcal G_r(\delta)$ is right $\C W_{r-2}$-linear, it induces a linear operator
\begin{equation}\label{eq:exact-induced-operator}
\mathcal G^{\bm\mu}= \mathcal G_r(\delta)\otimes_{\C W_{r-2}}\operatorname{Id}_{S^{\bm\mu}}:S^{(2,\bm\mu)}\longrightarrow S^{(2,\bm\mu)}.
\end{equation}
By Proposition~\ref{prop:auxiliary-operator-bimodule-map}, $\mathcal G^{\bm\mu}$ is a $\C W_r$-endomorphism. Let $\rho_{\bm\mu}:\C W_{r-2}\longrightarrow\operatorname{End}_{\C}(S^{\bm\mu})$ be the representation afforded by $S^{\bm\mu}$, then the matrix of
$\mathcal G^{\bm\mu}$ with respect to the basis $\{f(D_i)\otimes_{\C W_{r-2}} v_j\}$ is the block matrix
$$
\bigl(\rho_{\bm\mu}(\phi(D_i,D_j))\bigr)_{1\leq i,j\leq N}.
$$
Thus each entry of $\mathcal G_r(\delta)$, which
belongs to $\C W_{r-2}$, is replaced by the matrix representing its action on $S^{\bm\mu}$.

The operator $\mathcal G^{\bm\mu}$ must not be confused with the Gram matrix $G_{(2,\bm\mu)}$, which is defined using the bilinear form $\phi_{(2, \bm\mu)}$ on $S^{(2,\bm\mu)}$. Since $\mathcal G^{\bm\mu}$
is $\C W_r$-linear, it preserves each isotypic component of $S^{(2,\bm\mu)}$. On a multiplicity-one component isomorphic to $S^{\bm\lambda}$ it acts by a scalar, while on a multiplicity-two component it has the form $\operatorname{Id}_{S^{\bm\lambda}}\otimes A_{\bm\lambda,\bm\mu}(\delta)$, where $A_{\bm\lambda,\bm\mu}(\delta)$ is a $2\times2$ matrix acting on the vector space that records the two copies of $S^{\bm\lambda}$. We determine these actions in the next section.

\section{Homomorphisms into the second layer}

Throughout this section, we assume that $\delta\neq0$. Fix $\bm\mu\in\bm\Lambda_2(r-2)$ and let $\bm\lambda\in\bm\Lambda_2(r)$. We begin by introducing three $W_r$-invariant diagram sums. For $1\leq i\leq r$, let $E_i$ be the diagram with a horizontal edge joining $i$ and $-i$ in each row and straight edges elsewhere. For $1\leq i<j\leq r$, let $E_{ij}^{++}$ be the diagram with horizontal edges $\{i,j\}$ and $\{-i,-j\}$ in each row, and let $E_{ij}^{+-}$ be the diagram with horizontal edges $\{i,-j\}$ and $\{-i,j\}$ in each row. Set
$$
T_1=\sum_{i=1}^r E_i, \qquad T_2=\sum_{1\leq i<j\leq r}\bigl(E_{ij}^{++}+E_{ij}^{+-}\bigr),
$$
and
$$
T_3=\sum_{1\leq i<j\leq r}E_iE_j. 
$$
The three sums are invariant under conjugation by $W_r$. Hence their left actions are $W_r$-linear. The following lemma establishes relations that connect the calculations from our work \cite{KCTprevious} with the auxiliary operator constructed in Section~\ref{sec:fixed-bottom-operator}.

\begin{lemma}\label{lem:exact-second-layer-actions}
On $S^{(2,\bm\mu)}$, the following statements hold:
\begin{enumerate}
\item The operator $T_1$ maps $I_{01}\otimes_{\C W_{r-2}} S^{\bm\mu}$ into
$I_{20}\otimes_{\C W_{r-2}} S^{\bm\mu}$. Moreover, if
$\bm\lambda\in\mathcal A^+(\bm\mu)$, then the restriction
$$
T_1:M_{01}^{\bm\lambda,\bm\mu}
\longrightarrow M_{20}^{\bm\lambda,\bm\mu}
$$
is nonzero. If $\bm\lambda\in\mathcal A^-(\bm\mu)$, then $T_1M_{01}^{\bm\lambda,\bm\mu}=0$.
\item On $I_{01}\otimes_{\C W_{r-2}} S^{\bm\mu}$, the operator $T_2$ satisfies 
$$T_2=(\delta^2-2)\operatorname{Id}+\delta\tau+\Omega_r-\Omega_{r-2},$$
where $\tau$ is the involution defined in Lemma~\ref{lem:exact-tau-action}, $\Omega_r$ acts from the left, and $\Omega_{r-2}$ acts on the tensor factor $S^{\bm\mu}$;
\item The auxiliary operator $\mathcal G^{\bm\mu}$ satisfies 
$$\mathcal G^{\bm\mu}=T_3+T_2;$$
\item The operators $T_1$, $T_2$ and $T_3$ satisfy $$T_3=\frac12(T_1^2-\delta T_1).$$
\end{enumerate}
\end{lemma}

\begin{proof}
For (1), let $y\in D_{(0,1)}$. Let $y'\in D_{(2,0)}$ be obtained by replacing the two top horizontal edges
of $y$ by $\{i,-i\}$ and $\{j,-j\}$, leaving all other edges unchanged. By \cite[Lemma~5.2]{KCTprevious},
$$
T_1y=T_1\tau(y)=2y'.
$$
Consequently,
$$
T_1\bigl(y+\tau(y)\bigr)=4y',
\qquad
T_1\bigl(y-\tau(y)\bigr)=0.
$$
Every basis diagram in $D_{(2,0)}$ arises in this way. After tensoring with $S^{\bm\mu}$, these identities show that $T_1$ maps the $(+1)$-eigenspace of $\tau$ onto $I_{20}\otimes_{\C W_{r-2}}S^{\bm\mu}$ and annihilates the $(-1)$-eigenspace. By Lemma~\ref{lem:exact-tau-action}, the source and target of this surjective map each contain precisely
one copy of $S^{\bm\lambda}$ for every
$\bm\lambda\in\mathcal A^+(\bm\mu)$. Since $T_1$ is $\C W_r$-linear, Schur's lemma implies that its
restriction to $M_{01}^{\bm\lambda,\bm\mu}$ is nonzero. For $\bm\lambda\in\mathcal A^-(\bm\mu)$, the corresponding submodule lies in the $(-1)$-eigenspace and is therefore annihilated by $T_1$. This proves (1).

For~(2), see the the local calculation in \cite[Proposition~5.3]{KCTprevious}. The matching local configuration contributes $\delta^2$, while the other local configuration contributes $\delta\tau$. The terms obtained by moving one endpoint contribute $\Omega_r$. The terms involving only free vertices are counted once too, and their contribution is corrected by subtracting the right action of $\Omega_{r-2}$. The two repeated matching terms account for the constant $-2$.

For (3), let $D\in P(r,2)$ and set $X_D=D\otimes D\otimes1$. By the definition of $\phi$, for $E\in P(r,2)$ and $v\in S^{\bm\mu}$ we have
$$
X_D\bigl(f(E)\otimes_{\C W_{r-2}}v\bigr)=f(D)\otimes_{\C W_{r-2}}\phi(D,E)v.
$$
Consequently, as operators on $S^{(2,\bm\mu)}$,
$$
\mathcal G^{\bm\mu}=\sum_{D\in P(r,2)}X_D.
$$
If $D\in P_{(2,0)}$, then $X_D=E_iE_j$ for some
$i<j$. If $D\in P_{(0,1)}$, then $X_D$ is either
$E_{ij}^{++}$ or $E_{ij}^{+-}$. Hence
$$
\mathcal G^{\bm\mu}=\sum_{D\in P_{(2,0)}}X_D+\sum_{D\in P_{(0,1)}}X_D=T_3+T_2,
$$
which proves (3).

For (4), since $E_i^2=\delta E_i$ and
$E_iE_j=E_jE_i$ for $i\neq j$, we obtain
$$
T_1^2 =\sum_{i=1}^rE_i^2
+2\sum_{1\leq i<j\leq r}E_iE_j
=\delta T_1+2T_3.
$$
This proves (4).
\end{proof}

The decomposition of $I_r^{2'}(\delta)$ as a free right $\C W_{r-2}$-module gives the vector-space decomposition
\begin{equation}\label{eq:exact-cell-direct-sum}
S^{(2,\bm\mu)}\cong \bigoplus_{D\in P(r,2)} \left(f(D)\C W_{r-2}\otimes_{\C W_{r-2}}S^{\bm\mu}
\right).
\end{equation}
We shall use the following consequence of this decomposition in both admissible cases.

\begin{lemma}\label{lem:exact-top-layer-copy}
Let $M$ be a $\C W_r$-submodule of $S^{(2,\bm\mu)}$. If $T_1M=0$ and $\mathcal G^{\bm\mu}M=0$, then $J_1M=0$. In particular, $M$ is a $\B_r$-submodule, and if $M\cong S^{\bm\lambda}$ as a $\C W_r$-module, then $M\cong S^{(0,\bm\lambda)}$ as a $\B_r$-module.
\end{lemma}

\begin{proof}
By \eqref{eq:exact-cell-direct-sum}, every $m\in M$
can be written uniquely as
$$
m=\sum_{E\in P(r,2)}
f(E)\otimes_{\C W_{r-2}}v_E,
\qquad v_E\in S^{\bm\mu}.
$$
For $D\in P(r,2)$, recall that
$X_D=D\otimes D\otimes1$. Then
$$
X_D\cdot m
=
f(D)\otimes_{\C W_{r-2}}
\left(
\sum_{E\in P(r,2)}\phi(D,E)v_E
\right).
$$
Since $\mathcal G^{\bm\mu}=\sum_{D\in P(r,2)}X_D$
and $\mathcal G^{\bm\mu}m=0$, we obtain
$$
0=\sum_{D\in P(r,2)}X_D\cdot m.
$$
For distinct $D$, the vectors $X_D\cdot m$ lie in
distinct direct summands of
\eqref{eq:exact-cell-direct-sum}. Therefore, this
sum can vanish only if
$$
X_D\cdot m=0
\qquad\text{for every }D\in P(r,2).
$$
It follows that
$$
E_{ij}^{++}m=E_{ij}^{+-}m=E_iE_jm=0
\qquad(1\leq i<j\leq r).
$$
For a fixed $i$, the equality $T_1m=0$ implies
$$
0=E_iT_1m
=E_i^2m+\sum_{j\neq i}E_iE_jm
=\delta E_im.
$$
Since $\delta\neq0$, we have $E_iM=0$ for every $i$.

Thus the diagrams $E_i$, $E_{ij}^{++}$ and $E_{ij}^{+-}$ annihilate $M$. Together with $W_r$, these diagrams generate $\B_r$. Since $M$ is $\C W_r$-stable, it follows that $M$ is a $\B_r$-submodule.
The same diagrams generate $J_1$ as an ideal, so $J_1M=0$. Consequently, the action on $M$ factors through $\B_r/J_1\cong\C W_r$, which proves the final assertion.
\end{proof}

We now prove the converse of
Proposition~\ref{prop:known-local-restrictions}(b)
in the two cases corresponding to
$\mathcal A^-(\bm\mu)$ and $\mathcal A^+(\bm\mu)$.

\begin{theorem}\label{thm:exact-converse-minus}
Let $\bm\mu\in\bm\Lambda_2(r-2)$ and
$\bm\lambda\in\mathcal A^-(\bm\mu)$.
If $\delta^2-\delta-2+C_{\bm\lambda/\bm\mu}=0$, then
$$
\langle S^{(0,\bm\lambda)},
S^{(2,\bm\mu)}\rangle_r\neq0.
$$
\end{theorem}

\begin{proof}
Since $\bm\lambda\in\mathcal A^-(\bm\mu)$, the module
$I_{01}\otimes_{\C W_{r-2}}S^{\bm\mu}$ contains the unique copy
$$
M=M_{01}^{\bm\lambda,\bm\mu}\cong S^{\bm\lambda},
$$
whereas $I_{20}\otimes_{\C W_{r-2}}S^{\bm\mu}$ contains
no copy of $S^{\bm\lambda}$.
By Lemma~\ref{lem:exact-second-layer-actions}(1), the restriction
$$
T_1:M\longrightarrow I_{20}\otimes_{\C W_{r-2}}S^{\bm\mu}
$$
is a $\C W_r$-module homomorphism.
Since the target contains no constituent isomorphic to
$S^{\bm\lambda}$, it follows that $T_1M=0$.

By Lemma~\ref{lem:exact-tau-action}, $\tau$ acts as
$-\operatorname{Id}_M$.
Moreover, $\Omega_r$ acts on $M$ by the scalar
$C_{\bm\lambda}$, and $\Omega_{r-2}$ acts on
$S^{\bm\mu}$ by the scalar $C_{\bm\mu}$.
Therefore, Lemma~\ref{lem:exact-second-layer-actions}(2) gives
\begin{equation}\label{eq:exact-T2-minus}
T_2|_M=
\bigl(\delta^2-\delta-2+C_{\bm\lambda/\bm\mu}\bigr)
\operatorname{Id}_M.
\end{equation}
The hypothesis therefore implies that $T_2M=0$.
Since $T_1M=0$, Lemma~\ref{lem:exact-second-layer-actions}(4)
also yields $T_3M=0$.
Consequently,
$$
\mathcal G^{\bm\mu}M=(T_3+T_2)M=0.
$$
By Lemma~\ref{lem:exact-top-layer-copy},
$M\cong S^{(0,\bm\lambda)}$ as a $\B_r$-module.
The inclusion $M\hookrightarrow S^{(2,\bm\mu)}$ thus gives
the required nonzero homomorphism.
\end{proof}

In the case $\bm\lambda\in\mathcal A^+(\bm\mu)$,
the irreducible module $S^{\bm\lambda}$ occurs in both
$I_{20}\otimes_{\C W_{r-2}}S^{\bm\mu}$ and
$I_{01}\otimes_{\C W_{r-2}}S^{\bm\mu}$.
The proof requires a suitable linear combination of the embeddings into these two summands.

\begin{theorem}\label{thm:exact-converse-plus}
Let $\bm\mu\in\bm\Lambda_2(r-2)$ and
$\bm\lambda\in\mathcal A^+(\bm\mu)$.
If $\delta^2+\delta-4+C_{\bm\lambda/\bm\mu}=0$, then
$$
\langle S^{(0,\bm\lambda)},
S^{(2,\bm\mu)}\rangle_r\neq0.
$$
\end{theorem}

\begin{proof}
The $\bm\lambda$-isotypic component of $S^{(2,\bm\mu)}$
decomposes as
$M_{20}^{\bm\lambda,\bm\mu}\oplus
M_{01}^{\bm\lambda,\bm\mu}$,
where both summands are isomorphic to $S^{\bm\lambda}$.
Choose a nonzero $\C W_r$-module embedding
$$
\theta_{01}:S^{\bm\lambda}\longrightarrow
M_{01}^{\bm\lambda,\bm\mu}.
$$
Since $\bm\lambda\in\mathcal A^+(\bm\mu)$,
Lemma~\ref{lem:exact-tau-action} implies that $\tau$
acts as the identity on $M_{01}^{\bm\lambda,\bm\mu}$.
Moreover, by Lemma~\ref{lem:exact-second-layer-actions}(1),
the restriction of $T_1$ to
$M_{01}^{\bm\lambda,\bm\mu}$ is nonzero.
Define
$$
\theta_{20}=\frac12T_1\theta_{01}.
$$
Then $\theta_{20}$ is a nonzero $\C W_r$-module embedding into $M_{20}^{\bm\lambda,\bm\mu}$.
Every $\C W_r$-module homomorphism from $S^{\bm\lambda}$ into the $\lambda$-isotypic component can therefore be expressed as
$$
\Theta_{u,v}=u\theta_{20}+v\theta_{01}
$$
for some $u,v\in\C$.

By \cite[Lemma~5.2]{KCTprevious}, the action of $T_1$
on these embeddings is given by
$$
T_1\theta_{01}=2\theta_{20},
\qquad
T_1\theta_{20}=2\delta\theta_{20}.
$$
Consequently,
\begin{equation}\label{eq:exact-T1-plus}
T_1\Theta_{u,v}=2(\delta u+v)\theta_{20}.
\end{equation}
Thus, up to multiplication by a nonzero scalar,
the unique embedding whose image is annihilated by $T_1$ is
$$
\Theta=\theta_{20}-\delta\theta_{01}.
$$
Set $M=\operatorname{im}\Theta$.
By construction, $T_1M=0$.

By Lemma~\ref{lem:exact-tau-action}, $\tau$ acts as the identity on $M_{01}^{\bm\lambda,\bm\mu}$. Hence Lemma~\ref{lem:exact-second-layer-actions}(2) gives
$$
T_2\theta_{01}=
\bigl(\delta^2+\delta-2+C_{\bm\lambda/\bm\mu}\bigr)
\theta_{01}.
$$
Furthermore, with the above normalization of $\theta_{20}$,
the second formula in \cite[Proposition~5.3(ii)]{KCTprevious} gives
$$
T_2\theta_{20}=\delta(1+\tau)\theta_{01}
=2\delta\theta_{01}.
$$
Therefore,
\begin{equation}\label{eq:exact-T2-plus}
T_2\Theta
=-\delta
\bigl(\delta^2+\delta-4+C_{\bm\lambda/\bm\mu}\bigr)
\theta_{01}.
\end{equation}
The hypothesis implies that $T_2M=0$.
Since $T_1M=0$, Lemma~\ref{lem:exact-second-layer-actions}(4)
also gives $T_3M=0$.
Therefore,
$$
\mathcal G^{\bm\mu}M=(T_3+T_2)M=0.
$$

The preceding calculations may also be summarized by the
action matrix of the auxiliary operator.
By Lemma~\ref{lem:exact-second-layer-actions}(4), we have
$$
T_3\theta_{20}=\delta^2\theta_{20},
\qquad
T_3\theta_{01}=\delta\theta_{20}.
$$
Thus, with respect to the ordered pair
$(\theta_{20},\theta_{01})$, the action of
$\mathcal G^{\bm\mu}$ on the two copies of
$S^{\bm\lambda}$ is represented by the matrix
\begin{equation}\label{eq:exact-action-matrix}
A_{\bm\lambda,\bm\mu}(\delta)=
\begin{pmatrix}
\delta^2&\delta\\
2\delta&\delta^2+\delta-2+C_{\bm\lambda/\bm\mu}
\end{pmatrix}.
\end{equation}
More precisely, on the $\bm\lambda$-isotypic component,
we have
$$
\mathcal G^{\bm\mu}\big|_{
M_{20}^{\bm\lambda,\bm\mu}\oplus
M_{01}^{\bm\lambda,\bm\mu}}
=
\operatorname{Id}_{S^{\bm\lambda}}
\otimes A_{\bm\lambda,\bm\mu}(\delta).
$$
Under the identification of this component with \(S^{\bm\lambda}\) tensored with the two-dimensional multiplicity space, the embedding $\Theta=\theta_{20}-\delta\theta_{01}$
corresponds to the column vector
$(1,-\delta)^{\mathsf T}$. Direct multiplication gives
$$
A_{\bm\lambda,\bm\mu}(\delta)
\binom{1}{-\delta}
=
\binom{0}{
-\delta(\delta^2+\delta-4+C_{\bm\lambda/\bm\mu})}.
$$
The displayed matrix gives a second direct verification
that $\mathcal G^{\bm\mu}M=0$ under the stated quadratic
condition.
Together with $T_1M=0$,
Lemma~\ref{lem:exact-top-layer-copy} gives
$M\cong S^{(0,\bm\lambda)}$ as a $\B_r$-module.
Hence $\Theta$ defines the required nonzero
$\B_r$-homomorphism.
\end{proof}

Combining these converse results with Proposition~\ref{prop:known-local-restrictions}(b), we obtain a necessary and sufficient condition for the existence of a nonzero homomorphism from $S^{(0,\bm\lambda)}$ to $S^{(2,\bm\mu)}$.

\begin{corollary}\label{cor:exact-local-criterion}
Let $\bm\lambda\in\bm\Lambda_2(r)$ and
$\bm\mu\in\bm\Lambda_2(r-2)$.
Then:
\begin{enumerate}
\item If $\bm\lambda\in\mathcal A^-(\bm\mu)$, then
$$
\langle S^{(0,\bm\lambda)},S^{(2,\bm\mu)}\rangle_r\neq0
\quad\Longleftrightarrow\quad
\delta^2-\delta-2+C_{\bm\lambda/\bm\mu}=0.
$$
\item If $\bm\lambda\in\mathcal A^+(\bm\mu)$, then
$$
\langle S^{(0,\bm\lambda)},S^{(2,\bm\mu)}\rangle_r\neq0
\quad\Longleftrightarrow\quad
\delta^2+\delta-4+C_{\bm\lambda/\bm\mu}=0.
$$
\item If $\bm\lambda\notin\mathcal A(\bm\mu)$, then
$$
\langle S^{(0,\bm\lambda)},S^{(2,\bm\mu)}\rangle_r=0.
$$
\end{enumerate}
In cases~{\rm(1)} and~{\rm(2)}, whenever the Hom-space
is nonzero, it is one-dimensional.
\end{corollary}

\begin{proof}
The reverse implications in parts~(1) and~(2) follow
from Theorems~\ref{thm:exact-converse-minus} and
\ref{thm:exact-converse-plus}, respectively.
We prove the forward implications and establish the
remaining assertions.

Suppose that $\bm\lambda\in\mathcal A^-(\bm\mu)$
and that there exists a nonzero $\B_r$-module
homomorphism from $S^{(0,\bm\lambda)}$ to
$S^{(2,\bm\mu)}$.
Let $M$ be its image.
Since $M\cong S^{\bm\lambda}$ as a $\C W_r$-module
and $S^{\bm\lambda}$ occurs with multiplicity one
in $S^{(2,\bm\mu)}$, we must have
$M=M_{01}^{\bm\lambda,\bm\mu}$.
As $T_2\in J_1$ and $J_1M=0$, we have $T_2M=0$.
Equation~\eqref{eq:exact-T2-minus} therefore gives
$$
\delta^2-\delta-2+C_{\bm\lambda/\bm\mu}=0.
$$
The uniqueness of the copy of $S^{\bm\lambda}$,
together with Schur's lemma, implies that any two
nonzero homomorphisms differ by a nonzero scalar.
Thus, the Hom-space is one-dimensional.

Suppose that $\bm\lambda\in\mathcal A^+(\bm\mu)$
and that a nonzero $\B_r$-module homomorphism $\psi:S^{(0,\bm\lambda)}\longrightarrow S^{(2,\bm\mu)}$ exists. Identifying the source with $S^{\bm\lambda}$ as a
$\C W_r$-module, we may write
$\psi=u\theta_{20}+v\theta_{01}$, where $u,v\in\C$. Since the image of $\psi$ is annihilated by $J_1$,
it is annihilated by both $T_1$ and $T_2$. Equation~\eqref{eq:exact-T1-plus} then implies $v=-\delta u$.
Since $\psi\neq0$, we have $u\neq0$ and $\psi=u\Theta$. Applying Equation~\eqref{eq:exact-T2-plus}, we obtain
$$
0=-u\delta
\bigl(\delta^2+\delta-4+C_{\bm\lambda/\bm\mu}\bigr)
\theta_{01}.
$$
Since $u$, $\delta$, and $\theta_{01}$ are nonzero, we conclude that
$\delta^2+\delta-4+C_{\bm\lambda/\bm\mu}=0$. Moreover, the relation $v=-\delta u$ shows that every such homomorphism is a scalar multiple of $\Theta$. Thus, the Hom-space is one-dimensional.

Finally, suppose that
$\bm\lambda\notin\mathcal A(\bm\mu)$.
The restriction map
$$
\Hom_{\B_r}
\bigl(S^{(0,\bm\lambda)},S^{(2,\bm\mu)}\bigr)
\longrightarrow
\Hom_{\C W_r}
\bigl(S^{\bm\lambda},S^{(2,\bm\mu)}\bigr)
$$
is injective, and the second Hom-space is zero by the Proposition~\ref{prop:second-layer-Wr-decomposition}. This proves part~(3).
\end{proof}

\begin{remark}\label{rem:exact-matched-criterion}
The criterion depends on which component of the bipartition acquire the horizontal $2$-strip. If $\bm\lambda\in\mathcal A^-(\bm\mu)$, only the equation $\delta^2-\delta-2+C_{\bm\lambda/\bm\mu}=0$
is relevant for that pair. If $\bm\lambda\in\mathcal A^+(\bm\mu)$, the relevant equation is $\delta^2+\delta-4+C_{\bm\lambda/\bm\mu}=0$.
The criterion does not assert a converse for the other equation in either case.
\end{remark}

\section{Reduction to the second layer}

In this section, we prove that the semisimplicity of $\B_r$ is determined by the second-layer cell modules at rank $r$, where $r\geq 2$. Since the argument involves the tower of algebras, for $2\leq m\leq r$, we write $S_m^{(k,\bm\lambda)}$ for a cell module of $\B_m$ and $G_{(k,\bm\lambda)}^{\,m}$ for the Gram matrix of the bilinear form on $S_m^{(k,\bm\lambda)}$. Thus, $\det G_{(k,\bm\lambda)}^{\,m}=0$ if and only if $S_m^{(k,\bm\lambda)}$ is not simple. We denote its simple head by $D_m(k,\bm\lambda)$.

\subsection{Degeneracy of second-layer cell modules}
We begin by relating the degeneracy of a second-layer cell module to the homomorphisms studied in the preceding section.

\begin{lemma}
\label{lem:second-layer-degenerate-hom}
Let $\bm\mu\in\bm\Lambda_2(r-2)$. Then
$ \det G_{(2,\bm\mu)}^{\,r}=0$ if and only if there exists $\bm\lambda\in\bm\Lambda_2(r)$ such that
$
\left\langle S_r^{(0,\bm\lambda)},S_r^{(2,\bm\mu)}\right\rangle_r\neq0.
$
\end{lemma}

\begin{proof}
Suppose first that there is a nonzero homomorphism
$$
S_r^{(0,\bm\lambda)}
\longrightarrow S_r^{(2,\bm\mu)}.
$$
The source is simple, so its image is a simple submodule
of $S_r^{(2,\bm\mu)}$.
This image must be proper.
Otherwise, $S_r^{(2,\bm\mu)}$ would be simple and
isomorphic to $D_r(0,\bm\lambda)$, whereas its simple
head is $D_r(2,\bm\mu)$, and simple modules with
distinct cell labels are nonisomorphic.
Thus, $S_r^{(2,\bm\mu)}$ is not simple and its Gram
determinant is zero.

Conversely, suppose that $\det G_{(2,\bm\mu)}^{\,r}=0$.
Then the radical of $S_r^{(2,\bm\mu)}$ is nonzero.
Since $\B_r$ is finite-dimensional, this radical contains
a simple submodule, say $D_r(k,\bm\lambda)$.
Composing the canonical surjection
$$
S_r^{(k,\bm\lambda)}
\twoheadrightarrow D_r(k,\bm\lambda)
$$
with the inclusion into $S_r^{(2,\bm\mu)}$,
we obtain a nonzero homomorphism
$$
S_r^{(k,\bm\lambda)}
\longrightarrow S_r^{(2,\bm\mu)}.
$$
The unitriangular property of the decomposition matrix implies
$k\leq2$.

If $k=2$, then $\bm\lambda\neq\bm\mu$.
Indeed, the simple head $D_r(2,\bm\mu)$ occurs with multiplicity one in $S_r^{(2,\bm\mu)}$, by the unitriangular property of the decomposition matrix. Therefore, it cannot also occur in the radical.
By the rank-shifting property for homomorphisms between cell modules in Proposition~\ref{prop:shift-pairing}, we would then obtain
a nonzero homomorphism
$$
S_{r-2}^{(0,\bm\lambda)}
\longrightarrow S_{r-2}^{(0,\bm\mu)}
$$
between nonisomorphic simple modules, which is impossible.

If $k=1$, then again by Proposition~\ref{prop:shift-pairing}, we obtain
$$
\langle
S_{r-1}^{(0,\bm\lambda)},S_{r-1}^{(1,\bm\mu)}
\rangle_{r-1}\neq0,
$$
contrary to Proposition~\ref{prop:known-local-restrictions}(a),
since $\delta\neq0$.
Hence $k=0$, which proves the converse.
\end{proof}

Combining Lemma~\ref{lem:second-layer-degenerate-hom} with Corollary~\ref{cor:exact-local-criterion}, we obtain the following explicit form of the degeneracy condition.

\begin{corollary}
\label{cor:second-layer-degeneracy}
Let $\bm\mu\in\bm\Lambda_2(r-2)$. Then
$\det G_{(2,\bm\mu)}^{\,r}=0$
if and only if there exists $\bm\lambda\in\bm\Lambda_2(r)$ satisfying one of the following conditions:
\begin{enumerate}
\item $\bm\lambda\in\mathcal A^-(\bm\mu)$ and $\delta^2-\delta-2+C_{\bm\lambda/\bm\mu}=0$;
\item $\bm\lambda\in\mathcal A^+(\bm\mu)$ and $\delta^2+\delta-4+C_{\bm\lambda/\bm\mu}=0$.
\end{enumerate}
\end{corollary}

\subsection{Reduction through the tower}
We next establish a criterion for semisimplicity
in terms of the second-layer cell modules at all
ranks up to $r$.

\begin{proposition}\label{prop:exact-all-ranks}
The algebra $\B_r$ is semisimple if and only if
$\det G_{(2,\bm\mu)}^{\,m}\neq0$ for every
$2\leq m\leq r$ and every
$\bm\mu\in\bm\Lambda_2(m-2)$.
\end{proposition}

\begin{proof}
Suppose first that $\B_r$ is semisimple.
Assume, for a contradiction, that
$$
\det G_{(2,\bm\mu)}^{\,m}=0
$$
for some $2\leq m\leq r$ and
$\bm\mu\in\bm\Lambda_2(m-2)$.
If $m=r$, this contradicts the semisimplicity of $\B_r$,
so we may assume that $m<r$.
By Corollary~\ref{cor:second-layer-degeneracy},
there exists $\bm\lambda\in\bm\Lambda_2(m)$ satisfying one of the two matched conditions stated there. Set $\ell=r-m$. For a partition $\alpha=(\alpha_1,\alpha_2,\ldots)$, define
$$
\alpha^{[\ell]}=(\alpha_1+\ell,\alpha_2,\ldots),
$$
with $\varnothing^{[\ell]}=(\ell)$.

Suppose first that $\bm\lambda\in\mathcal A^-(\bm\mu)$. Write $\alpha=\lambda^{(1)}=\mu^{(1)}$, and set
$$
\widetilde{\bm\lambda}
=\bigl(\alpha^{[\ell]},\lambda^{(2)}\bigr),
\qquad
\widetilde{\bm\mu}
=\bigl(\alpha^{[\ell]},\mu^{(2)}\bigr).
$$
Then $\widetilde{\bm\lambda}\in\bm\Lambda_2(r)$,
$\widetilde{\bm\mu}\in\bm\Lambda_2(r-2)$, and
$\widetilde{\bm\lambda}\in
\mathcal A^-(\widetilde{\bm\mu})$.
Since the skew diagram lies entirely in the unchanged
second component, we have
$$
C_{\widetilde{\bm\lambda}/\widetilde{\bm\mu}}
=C_{\bm\lambda/\bm\mu}.
$$
Thus, the same quadratic equation holds for the lifted pair.
By Corollary~\ref{cor:second-layer-degeneracy},
we have
$$
\det G_{(2,\widetilde{\bm\mu})}^{\,r}=0,
$$
contradicting the semisimplicity of $\B_r$.

Now suppose that $\bm\lambda\in\mathcal A^+(\bm\mu)$.
Write $\beta=\lambda^{(2)}=\mu^{(2)}$ and define
$$
\widetilde{\bm\lambda}
=\bigl(\lambda^{(1)},\beta^{[\ell]}\bigr),
\qquad
\widetilde{\bm\mu}
=\bigl(\mu^{(1)},\beta^{[\ell]}\bigr).
$$
Now the skew diagram lies entirely in the unchanged
first component.
Thus,
$\widetilde{\bm\lambda}\in
\mathcal A^+(\widetilde{\bm\mu})$ and
$$
C_{\widetilde{\bm\lambda}/\widetilde{\bm\mu}}
=C_{\bm\lambda/\bm\mu}.
$$
Therefore, the corresponding quadratic equation is
preserved under this lifting as well.
Corollary~\ref{cor:second-layer-degeneracy}
again gives a singular second-layer Gram matrix at
rank $r$, contradicting the semisimplicity of $\B_r$.
This proves the forward implication.

Conversely, assume that all second-layer Gram determinants
at ranks $2,\ldots,r$ are nonzero, but that $\B_r$
is not semisimple.
Then there exists a cell module $S_r^{(k_2,\bm\mu)}$
with a nonzero radical.
As in the proof of
Lemma~\ref{lem:second-layer-degenerate-hom},
a simple submodule of this radical yields a cell label
$(k_1,\bm\lambda)\neq(k_2,\bm\mu)$ such that
$$
\langle
S_r^{(k_1,\bm\lambda)},S_r^{(k_2,\bm\mu)}
\rangle_r\neq0.
$$
The partial order in the cell datum implies that
$k_1\leq k_2$.
If $k_1=k_2$, then $\bm\lambda\neq\bm\mu$.
Using Proposition~\ref{prop:shift-pairing} would then give
a nonzero homomorphism between the nonisomorphic
simple modules $S_{r-k_1}^{(0,\bm\lambda)}$ and
$S_{r-k_1}^{(0,\bm\mu)}$, which is impossible.
Hence $k_1<k_2$.

Set $k=k_2-k_1$ and $n=r-k_1$.
Applying Proposition~\ref{prop:shift-pairing} once more, we obtain
$$
\langle
S_n^{(0,\bm\lambda)},S_n^{(k,\bm\mu)}
\rangle_n\neq0.
$$
Starting from this homomorphism, we repeatedly apply
the reduction to smaller rank in Proposition~\ref{prop:case-split}, retaining at each step
only nonzero homomorphisms whose targets have positive
layer indices.
Among these homomorphisms, including the original one,
choose one of minimal rank, say
$$
\langle
S_m^{(0,\bm\alpha)},S_m^{(q,\bm\beta)}
\rangle_m\neq0
$$
with $q>0$.
If $q\geq3$, Proposition~\ref{prop:case-split} gives
a nonzero homomorphism at rank $m-1$ whose target
lies in layer $q$, $q-1$, or $q-2$.
All three layer indices are positive, contradicting
the minimality of $m$.
Therefore, $q=1$ or $q=2$.

The case $q=1$ contradicts Proposition~\ref{prop:known-local-restrictions}(a).
Hence $q=2$, and
Lemma~\ref{lem:second-layer-degenerate-hom} implies that
$$
\det G_{(2,\bm\beta)}^{\,m}=0,
$$
contrary to our assumption.
Thus $\B_r$ is semisimple.
\end{proof}

The lifting construction used to prove the forward
implication of Proposition~\ref{prop:exact-all-ranks}
yields the following consequence.

\begin{corollary}\label{cor:exact-lift-degeneracy}
If $\det G_{(2,\bm\mu)}^{\,m}=0$ for some
$2\leq m\leq r$ and $\bm\mu\in\bm\Lambda_2(m-2)$,
then there exists
$\widetilde{\bm\mu}\in\bm\Lambda_2(r-2)$ such that
$\det G_{(2,\widetilde{\bm\mu})}^{\,r}=0$.
\end{corollary}

\begin{proof}
If $m=r$, take $\widetilde{\bm\mu}=\bm\mu$.
If $m<r$, the required bipartition is obtained from
the lifting construction used to prove the forward
implication of Proposition~\ref{prop:exact-all-ranks}.
\end{proof}

We now state the semisimplicity criterion in terms of
the second-layer cell modules at rank $r$.

\begin{theorem}\label{thm:exact-deciding-layer}
Let $r\geq2$ and $\delta\neq0$.
Then $\B_r$ is semisimple if and only if
$\det G_{(2,\bm\mu)}^{\,r}\neq0$ for every
$\bm\mu\in\bm\Lambda_2(r-2)$.
\end{theorem}

\begin{proof}
If $\B_r$ is semisimple, all its cell modules are
simple, so the stated determinants are nonzero.
Conversely, suppose that all second-layer Gram
determinants at rank $r$ are nonzero.
If $\B_r$ were not semisimple, then, by
Proposition~\ref{prop:exact-all-ranks}, there would
exist an integer $m$ with $2\leq m\leq r$ and
a bipartition $\bm\mu\in\bm\Lambda_2(m-2)$ such that
$\det G_{(2,\bm\mu)}^{\,m}=0$.
Corollary~\ref{cor:exact-lift-degeneracy} would then
give a singular second-layer Gram matrix at rank $r$,
a contradiction.
\end{proof}

\section{The semisimplicity criterion and the exceptional parameters}
\label{sec:semisimplicity-criterion}

We now combine the homomorphism criterion obtained in Corollary~\ref{cor:exact-local-criterion} with Theorem~\ref{thm:exact-deciding-layer} from the preceding section. It remains to determine the possible values of $C_{\bm\lambda/\bm\mu}$ for $\bm\mu\in\bm\Lambda_2(r-2)$ and
$\bm\lambda\in\mathcal A(\bm\mu)$.

\subsection{Possible content sums of horizontal $2$-strips}

\begin{proposition}\label{prop:exact-content-spectrum}
For $r\geq2$, define
\begin{equation*}
\mathcal C_r^-=
\left\{
C_{\bm\lambda/\bm\mu}:
\bm\mu\in\bm\Lambda_2(r-2),\
\bm\lambda\in\mathcal A^-(\bm\mu)
\right\}
\,\text{ and }\,
\mathcal C_r^+=
\left\{
C_{\bm\lambda/\bm\mu}:
\bm\mu\in\bm\Lambda_2(r-2),\
\bm\lambda\in\mathcal A^+(\bm\mu)
\right\}.
\end{equation*}
Then $\mathcal C_r^-=\mathcal C_r^+=
\{2s:s\in\mathscr S_r\}$, where $\mathscr S_r$
is the set defined in \eqref{eq:exact-index-set}.
\end{proposition}

\begin{proof}
Recall that
$$
C_{\bm\lambda/\bm\mu}
=2\sum_{p\in[\bm\lambda/\bm\mu]}c(p).
$$
We first determine the possible values of the unscaled sum
$$
s_{\bm\lambda/\bm\mu}
=\sum_{p\in[\bm\lambda/\bm\mu]}c(p).
$$

Suppose that $\bm\lambda\in\mathcal A^-(\bm\mu)$.
Then $\lambda^{(1)}=\mu^{(1)}$ and
$\lambda^{(2)}/\mu^{(2)}$ is a horizontal $2$-strip.
Set $\alpha=\lambda^{(2)}$, $\beta=\mu^{(2)}$,
and $n=|\alpha|$.
Then $2\leq n\leq r$, $|\beta|=n-2$, and
$\alpha/\beta$ is an ordinary horizontal $2$-strip.
Thus, every content sum arising from a pair
$(\bm\lambda,\bm\mu)$ with
$\bm\lambda\in\mathcal A^-(\bm\mu)$ is the content sum
of an ordinary horizontal $2$-strip $\alpha/\beta$
with $2\leq|\alpha|\leq r$.

Conversely, let $\alpha/\beta$ be an ordinary
horizontal $2$-strip with $|\alpha|=n$,
$|\beta|=n-2$, and $2\leq n\leq r$.
Define the partition $\gamma$ by $\gamma=(r-n)$
if $n<r$, and $\gamma=\varnothing$ if $n=r$.
Now set
$$
\bm\lambda=(\gamma,\alpha),
\qquad
\bm\mu=(\gamma,\beta).
$$
Then $\bm\lambda\in\mathcal A^-(\bm\mu)$ and the
content sum is unchanged.
Consequently, the content sums arising from pairs
with $\bm\lambda\in\mathcal A^-(\bm\mu)$ are precisely
the content sums of ordinary horizontal $2$-strips
$\alpha/\beta$ satisfying $2\leq|\alpha|\leq r$.

For a bipartition
$\bm\nu=(\nu^{(1)},\nu^{(2)})$, set
$\bm\nu^{\mathrm{op}}=(\nu^{(2)},\nu^{(1)})$.
The map
$$
(\bm\lambda,\bm\mu)
\longmapsto
(\bm\lambda^{\mathrm{op}},\bm\mu^{\mathrm{op}})
$$
is a bijection between the admissible pairs satisfying
$\bm\lambda\in\mathcal A^-(\bm\mu)$ and those satisfying
$\bm\lambda\in\mathcal A^+(\bm\mu)$.
Since interchanging the two components preserves
the contents of the added boxes, the two sets of
content values coincide.

It remains to determine the possible sums of contents
of ordinary horizontal $2$-strips.
Following the argument in
\cite[proof of Proposition~4.9(1)]{RX07}, define
$$
Z(r)=
\left\{
1-\sum_{p\in[\alpha/\beta]}c(p)
\;\middle|\;
\begin{array}{c}
|\alpha|=n,\ |\beta|=n-2,\ 2\leq n\leq r,\\
\alpha/\beta\text{ is a horizontal $2$-strip}
\end{array}
\right\}.
$$
By \cite[Proposition~2.4]{RS06},
$$
Z(r)=
\{i\in\Z:4-2r\leq i\leq r-2\}
\setminus
\{i\in\Z:4-2r<i\leq3-r,\ 2\nmid i\}.
$$
Set $s=1-i$.
The inequalities $4-2r\leq i\leq r-2$ then become
$3-r\leq s\leq2r-3$.
Moreover, the excluded odd values of $i$ correspond
exactly to the even values of $s$ satisfying
$r-2\leq s<2r-3$.
Hence the set of possible values of $s$ is
$$
\{s\in\Z:3-r\leq s\leq r-3\}
\cup\{1,3,\ldots,2r-3\}
=\mathscr S_r.
$$
Since $C_{\bm\lambda/\bm\mu}
=2s_{\bm\lambda/\bm\mu}$, multiplying these possible
values of $s_{\bm\lambda/\bm\mu}$ by $2$
proves the result.
\end{proof}

\subsection{Proof of Theorem~A}

\begin{proof}
Suppose first that $\B_r$ is semisimple.
We show that condition~{\rm(2)} holds.
Assume, to the contrary, that it fails.
Then there exists $s\in\mathscr S_r$ for which
$\delta^2-\delta-2+2s=0$ or
$\delta^2+\delta-4+2s=0$.
In the first case,
Proposition~\ref{prop:exact-content-spectrum}
guarantees the existence of
$\bm\mu\in\bm\Lambda_2(r-2)$ and
$\bm\lambda\in\mathcal A^-(\bm\mu)$ with
$C_{\bm\lambda/\bm\mu}=2s$.
In the second case, it gives
$\bm\mu\in\bm\Lambda_2(r-2)$ and
$\bm\lambda\in\mathcal A^+(\bm\mu)$ with
$C_{\bm\lambda/\bm\mu}=2s$.
In either case, Corollary~\ref{cor:second-layer-degeneracy}
implies that
$$
\det G_{(2,\bm\mu)}^{\,r}=0,
$$
contrary to Theorem~\ref{thm:exact-deciding-layer}.
Therefore, condition~{\rm(2)} is necessary.

Conversely, assume that condition~{\rm(2)} holds.
Suppose, for a contradiction, that $\B_r$ is not
semisimple.
By Theorem~\ref{thm:exact-deciding-layer}, there exists
$\bm\mu\in\bm\Lambda_2(r-2)$ such that
$\det G_{(2,\bm\mu)}^{\,r}=0$.
By Corollary~\ref{cor:second-layer-degeneracy},
there exists $\bm\lambda\in\bm\Lambda_2(r)$ such that
either $\bm\lambda\in\mathcal A^-(\bm\mu)$ and
$\delta^2-\delta-2+C_{\bm\lambda/\bm\mu}=0$,
or $\bm\lambda\in\mathcal A^+(\bm\mu)$ and
$\delta^2+\delta-4+C_{\bm\lambda/\bm\mu}=0$.
By Proposition~\ref{prop:exact-content-spectrum},
$C_{\bm\lambda/\bm\mu}=2s$ for some
$s\in\mathscr S_r$.
Hence one of the two expressions in
condition~{\rm(2)} must vanish, a contradiction.
Therefore, $\B_r$ is semisimple.
\end{proof}

\subsection{The cases $r=2$ and $r=3$}

We now illustrate the criterion in
Theorem~A for the first two
nontrivial cases.

For $r=2$, we have $\mathscr S_2=\{1\}$.
The two expressions in Theorem~A are
$$
\delta^2-\delta=\delta(\delta-1),
\qquad
\delta^2+\delta-2=(\delta-1)(\delta+2).
$$
Since $\delta\neq0$,
$$
\B_2(C,\delta)\text{ is semisimple}
\quad\Longleftrightarrow\quad
\delta\notin\{1,-2\}.
$$

For $r=3$, we have $\mathscr S_3=\{0,1,3\}$.
The corresponding polynomials are listed below.
$$
\begin{array}{c|c|c}
s&\delta^2-\delta-2+2s&\delta^2+\delta-4+2s\\ \hline
0&\delta^2-\delta-2&\delta^2+\delta-4\\
1&\delta^2-\delta&\delta^2+\delta-2\\
3&\delta^2-\delta+4&\delta^2+\delta+2
\end{array}
$$
Solving the corresponding equations and excluding
$\delta=0$ gives
$$
\begin{aligned}
\B_3(C,\delta)\text{ is not semisimple}
\quad\Longleftrightarrow\quad
\delta\in\Bigg\{
2,-1,1,-2,
\frac{-1\pm\sqrt{17}}{2},
\frac{1\pm\mathrm i\sqrt{15}}{2},
\frac{-1\pm\mathrm i\sqrt7}{2}
\Bigg\}.
\end{aligned}
$$

\begin{remark}
Throughout this paper we work over $\C$ and assume that $\delta\neq0$. The positive-characteristic case requires separate treatment. Let $K$ be a field of characteristic $p>0$ and let $\delta\in K^\times$. If $p\mid |W_r|=2^r r!$, then Maschke's theorem implies that $K W_r$ is not semisimple. Since 
$$
\B_r(C, \delta)/J_1\cong KW_r,
$$ 
the algebra $\B_r(C,\delta)$ cannot be semisimple in this case. Even
when $p\nmid 2^r r!$, the branching rules, the decomposition of the
second-layer cell modules, and the homomorphism calculations used above require separate verification in positive characteristic. We therefore restrict the present paper to the complex field.  
\end{remark}

\end{document}